\documentclass[11pt,a4paper,oneside,reqno]{amsart}
\usepackage{lmodern}
\usepackage[T1]{fontenc}
\usepackage[english]{babel}
\usepackage{amsopn}
\usepackage{graphicx} 
\usepackage{amssymb}
\usepackage{booktabs}
\usepackage{siunitx} 
\usepackage{mathrsfs}  

\usepackage{hyperref}

\newtheorem{theorem}{Theorem}[section]

\newtheorem{question}[theorem]{Question}

\newtheorem{conjecture}[theorem]{Conjecture}

\newtheorem{remark}[theorem]{Remark}

\title[FAMED by computer: the A.K. volume conjecture for 42,000 knots]{
FAMED by computer: proving the Andersen-Kashaev volume conjecture for 42,000 knots}
\author{Fathi Ben Aribi, Antonin Guilloux and Ka Ho Wong}

\newcommand{\C}{\mathbb{C}}
\newcommand{\R}{\mathbb{R}}

\newcommand{\Z}{\mathbb{Z}}

\newcommand{\Log}{\mathrm{Log}}
\newcommand{\Vol}{\mathrm{Vol}}

\usepackage{brunnian} 
\usetikzlibrary{knots}

\usetikzlibrary{matrix}
\usetikzlibrary{arrows}

\usetikzlibrary{patterns}
\usetikzlibrary{shapes}
\usetikzlibrary{arrows}
\usetikzlibrary{decorations.markings}
\usetikzlibrary{decorations.pathreplacing}

\tikzset{
  on each segment/.style={
    decorate,
    decoration={
      show path construction,
      moveto code={},
      lineto code={
        \path [#1]
        (\tikzinputsegmentfirst) -- (\tikzinputsegmentlast);
      },
      curveto code={
        \path [#1] (\tikzinputsegmentfirst)
        .. controls
        (\tikzinputsegmentsupporta) and (\tikzinputsegmentsupportb)
        ..
        (\tikzinputsegmentlast);
      },
      closepath code={
        \path [#1]
        (\tikzinputsegmentfirst) -- (\tikzinputsegmentlast);
      },
    },
  },
  mid arrow/.style={postaction={decorate,decoration={
        markings,
        mark=at position .5 with {\arrow[#1]{>}}
      }}}
      ,
}

\tikzset{
  on each segment/.style={
    decorate,
    decoration={
      show path construction,
      moveto code={},
      lineto code={
        \path [#1]
        (\tikzinputsegmentfirst) -- (\tikzinputsegmentlast);
      },
      curveto code={
        \path [#1] (\tikzinputsegmentfirst)
        .. controls
        (\tikzinputsegmentsupporta) and (\tikzinputsegmentsupportb)
        ..
        (\tikzinputsegmentlast);
      },
      closepath code={
        \path [#1]
        (\tikzinputsegmentfirst) -- (\tikzinputsegmentlast);
      },
    },
  },
  mid arrow d/.style={postaction={decorate,decoration={
        markings,
        mark=at position .5 with {\arrow[#1]{>>}}
      }}}
      ,
}

\tikzset{
  on each segment/.style={
    decorate,
    decoration={
      show path construction,
      moveto code={},
      lineto code={
        \path [#1]
        (\tikzinputsegmentfirst) -- (\tikzinputsegmentlast);
      },
      curveto code={
        \path [#1] (\tikzinputsegmentfirst)
        .. controls
        (\tikzinputsegmentsupporta) and (\tikzinputsegmentsupportb)
        ..
        (\tikzinputsegmentlast);
      },
      closepath code={
        \path [#1]
        (\tikzinputsegmentfirst) -- (\tikzinputsegmentlast);
      },
    },
  },
  mid arrow s/.style={postaction={decorate,decoration={
        markings,
        mark=at position .5 with {\arrow[#1]{stealth}}
      }}}
      ,
}

\tikzset{
  on each segment/.style={
    decorate,
    decoration={
      show path construction,
      moveto code={},
      lineto code={
        \path [#1]
        (\tikzinputsegmentfirst) -- (\tikzinputsegmentlast);
      },
      curveto code={
        \path [#1] (\tikzinputsegmentfirst)
        .. controls
        (\tikzinputsegmentsupporta) and (\tikzinputsegmentsupportb)
        ..
        (\tikzinputsegmentlast);
      },
      closepath code={
        \path [#1]
        (\tikzinputsegmentfirst) -- (\tikzinputsegmentlast);
      },
    },
  },
  mid arrow l/.style={postaction={decorate,decoration={
        markings,
        mark=at position .5 with {\arrow[#1]{latex}}
      }}}
      ,
}

\usepackage{graphicx}

\usepackage{nicematrix}

\newcommand{\fathi}[1]{\textcolor{magenta}{#1}}

\begin{document}

\begin{abstract}
    The FAMED condition is a combinatorial property for ideal triangulations of $3$-manifolds, which was introduced in 2024 by the first and last authors in order to study the Andersen--Kashaev volume conjecture. They notably proved  that this conjecture is true for all FAMED geometric triangulations of one-cusped hyperbolic $3$-manifolds with trivial second homology.

    In this paper, using a straightforward computer implementation in \textit{Regina} and \textit{Snappy}, we find FAMED geometric triangulations for more than 42~000  complements of knots in $S^3$, including all knots with $12$ crossings or fewer and all knots whose complement can be triangulated with $23$ tetrahedra or fewer. 
    As a consequence, the Andersen--Kashaev conjecture is now proven to be true for as many new examples.    

    Along the way, we find several new insights about the FAMED property, which have great value in the quest of a general proof of the Andersen--Kashaev volume conjecture for every knot complement.
\end{abstract}

\maketitle

\section{Introduction}

For any one-cusped hyperbolic $3$-manifold $M$ (such as the complement $M=S^3 \setminus K$ of a hyperbolic knot $K$), a quantum invariant $(Z_\hbar(M))_{\hbar>0}$ is a certain kind of numerical topological invariant of $M$, often made from the meeting of a combinatorial or algebraic construction and a physical inspiration. 

The associated \textit{volume conjecture} for this invariant surmises that $Z_\hbar(M)$ should manifest an exponential growth in the semi-classical limit $\hbar \to 0^+$, where the growth rate should be the hyperbolic volume $\mathrm{Vol}(M)$ of the manifold. We refer to \cite{MY} for a survey on volume conjectures.

Volume conjectures are beautifully diverse in their statements and strategies of proofs, as they encompass fields such as combinatorics, topology, hyperbolic geometry, quantum algebra and complex asymptotical analysis. Notably, proving a volume conjecture usually requires technical analytical tools such as the saddle point method. The technical conditions needed to apply this method can usually be verified only in certain specific examples.

Usually.
But \textit{not always.} Not anymore.
Indeed, in \cite{BAW}, for the Andersen-Kashaev volume conjecture, the first and last authors proved that the analytical convergence properties directly follow from specific geometric and combinatorial properties. Hence the proof of the Andersen--Kashaev volume conjecture reduces to checking these geometric and combinatorial properties.

Let us clarify. It is often useful to define a quantum invariant $Z_\hbar(M)$
from an ordered ideal triangulation $X$ of the 3-manifold $M$, usually as a state sum of functions on the tetrahedra of $X$, where the functions contain quantum deformations of the classical dilogarithm 
$$ \mathrm{Li}_2(z) := - \int_0^z \Log(1-u) \frac{du}{u} \ \ \ \textrm{for} \ z \in \C \setminus [1,\infty).$$
As the volume functional associated to the triangulation $X$ is constructed for such dilogarithms, one can use asymptotical analysis to prove that when $\hbar\to0^+$, the invariant $Z_\hbar(M)$ grows or decreases as $\exp\left(\frac{\pm 1}{2\pi \hbar}\Vol(M)\right)$.

In \cite{BAW}, the first and last authors introduced the \textit{FAMED condition}, which is a combinatorial property for an ordered ideal triangulation $X$. They then proved that for the Andersen--Kashaev TQFT invariant $Z_\hbar(M)$, the associated Andersen--Kashaev volume conjecture is true for the triangulation $X$ of the manifold $M$ if $X$ is FAMED and geometric (where \textit{geometric} means that $X$ admits a positive angle structure that encodes the unique complete hyperbolic structure on $M$). Moreover, checking if a given triangulation $X$ is FAMED is easy: one has to compute straightforward matrices associated to faces and edges and check that they satisfy an explicit algebraic relation (see Section \ref{sub:FAMED}).

Thus, the reader can now be either disappointed or relieved to learn that we will not mention more quantum topology nor complex analysis in the remainder of this paper. The goal of this paper is to detail how one can check the FAMED condition on the computer and to present the results we obtained on the Hoste--Thistlethwaite census of triangulations of knot complements, as integrated into \textit{\textit{Snappy}} under the name \textit{HTLinks}. These results can be summarized as follows:

\begin{theorem}\label{thm:main}
\begin{enumerate}
    \item Among the 59~937 triangulations of knot complements in the \textit{HTlinks} \textit{\textit{Snappy}} census, for 42~117 of these examples, there is a triangulation of the same knot complement that is FAMED and certified geometric.
    
    These examples include all hyperbolic knots with 12 crossings or fewer and all but 6 knots whose initial \textit{HTLinks} triangulation has 23 tetrahedra or fewer.
    \item In all examples, the triangulation is FAMED for the preferred longitude $l_K$ of the knot $K$.
    \item \textbf{For these 42~117 triangulations of distinct knot complements, the Andersen--Kashaev volume conjecture is proven.}
\end{enumerate}
    \end{theorem}

    {Theorem \ref{thm:main} (3) is the main result of this paper, and illustrates that the FAMED condition holds promise for the study of volume conjectures.}

The process of this computation is straightforward: for each example, we use built-in \textit{Regina} \cite{Regina} and \textit{\textit{Snappy}} \cite{Snappy} methods or custom functions to:
\begin{enumerate}
    \item Try to order (using \textit{Regina}) or find all orders for the given triangulation.
    \item Try and certify that the triangulation is geometric, using \textit{Snappy}'s certified ball arithmetic methods.
    \item For the orders found in Step 1, check the combinatorial FAMED condition by computing the set of matrices described in Section \ref{sec:prelim}. This is done via a mix of built-in \textit{Snappy} and \textit{Regina} methods and some custom methods tailored for this specific goal.
    \item If one FAMED and geometric triangulation is found, the task is done;
    \item Else, use \textit{Regina}'s retriangulate method to explore all retriangulations of the initial triangulation (within a given distance in the graph of retriangulation) to try and find a FAMED and geometric triangulation of the same manifold.
\end{enumerate}
We chose to perform these computations using SageMath \cite{sagemath} with \textit{Snappy} and \textit{Regina} added and wrote a Sage module \cite{FAMEDexploration}. Depending on the strategy, in the range of complexity of \textit{Snappy}'s \textit{HTLinks} Census, the first four steps may take up to a few seconds, when trying to list all possible orders of a triangulation. The last step instead is not bounded in complexity as is and may take a very long time and memory space. At this point, we lack a rationale to guide our retriangulation process.

We were able to perform the computation for all manifolds in the census up to 23 tetrahedra, using parallel computing. We arbitrarily stopped trying to retriangulate at this point because the computing time until we found a FAMED geometric triangulation grew seemingly doubly exponentially as the initial number of tetrahedra in the \textit{HTLinks} triangulation grew (also, $42$ is a pretty number). However, we have little doubt that given enough computing power and time - or a better retriangulate strategy, one could find a positive answer to the following question for the whole \textit{HTLinks} census:

\begin{question}
\begin{enumerate}
    \item Does every one-cusped hyperbolic $3$-manifold admit a FAMED geometric triangulation?
    \item (Weaker question:) Does every hyperbolic knot complement in $\mathbb{S}^3$ admit a FAMED geometric triangulation?
\end{enumerate}
\end{question}

These experimental (but nevertheless exact) results lead us to believe that FAMED triangulations have great potential for studying volume conjectures, as long as similar results as those of \cite{BAW} can be established for other quantum invariants.

In Section \ref{sec:prelim}, we review necessary preliminaries such as the FAMED condition. In Section \ref{sec:coding}, we detail the steps of our coding process. In Section \ref{sec:patterns} we list some interesting phenomena we observed in the process and the final results. Finally, in Section \ref{sec:future} we review future prospects and generalizations of our techniques.

\section*{Acknowledgements}

The authors thank François Guéritaud, Clément Maria, Stéphane Baseilhac and Ryan Budney for stimulating discussions. 

They also thank the communities of developers for Regina, SageMath, and Snappy, as well as IMJ-PRG and Ouragan INRIA for providing computing capabilities and servers.

\section{Mathematical preliminaries}\label{sec:prelim}

We will mostly follow the notation of \cite{BAW}.
In the rest of this paper, $M$ will denote a one-cusped hyperbolic $3$-manifold with trivial second homology. Two remarks are in order:
\begin{itemize}
    \item Firstly, most of the results we present here (Theorem \ref{thm:main}) concern complements  $M=S^3\setminus K$ of hyperbolic knots $K$ in $S^3$, so the reader can think $M$ as a knot complement for most of the paper without harm.
    \item Secondly, the restrictions on the number of cusps and the second homology come from technical requirements to define the Andersen--Kashaev TQFT \cite{AK}, but could be generalized in future works.
\end{itemize}

Likewise, $X$ will denote an ideal triangulation with $N$ tetrahedra $T_1$,$ \ldots$, $T_N$ and a gluing relation $\sim$ that leaves none of the $4N$ initial faces unpaired. The associated sets of cells will be denoted
\begin{itemize}
    \item $X^3 = \{T_1, \ldots, T_N\}$ the set of $3$-cells made of the interiors of the tetrahedra $T_i$ (which will also be denoted $T_i$ when there is no risk of confusion),
    \item $X^2$ the set of $2N$ quotient faces,
    \item $X^1$ the set of $N$ quotient edges,
    \item $X^0=\{\star\}$ the single quotient vertex, which represents the collapsed knot $K$ when $M=S^3 \setminus K$.
\end{itemize}

\begin{figure}[!h]
    \centering
    \includegraphics[scale=.5]{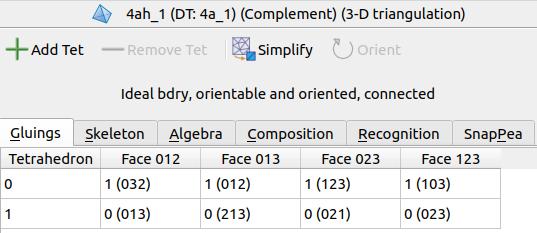}
    \caption{A \textit{\textit{Regina}} triangulation of the complement of the figure-eight knot $4_1$}
    \label{fig:Regina:41}
\end{figure}

We mostly use \textit{Regina} software to study triangulations. An ideal triangulation $X$ is displayed in \textit{Regina} as in Figure \ref{fig:Regina:41} (where $M=S^3\setminus 4_1$ and $N=2$). The tetrahedra are numbered from $0$ to $N-1$, and for each one, its $4$ (pre-gluing) vertices are numbered from $0$ to $3$, its $6$ edges are named $01, 02,03,12,13, 23$, and its $4$ faces are named $012,013,023,123$. 
Let us illustrate the way gluing is encoded with an example. In Figure \ref{fig:Regina:41}, the face $012$ of tetrahedron $0$ is glued to the face $023$ of the tetrahedron $1$, but without respecting the order of vertices: $0$ is sent to $0$, $1$ to $3$ and $2$ to $2$. This defines a unique affine gluing of both triangular faces, encoded by $1(032)$ and $0(021)$ in the corresponding table cells.

\subsection{Ordered triangulations}\label{sub:order}

We now say that a triangulation $X$ is \textit{ordered} when
it is endowed with a vertex numbering $0,1,2,3$ on each of its tetrahedra so that
the face gluing respects the order of the vertices (the smallest one of a  is always sent to the smallest one, and the highest to the highest). We say $X$ is \textit{orderable} when it can be ordered.

\begin{figure}[!h]
    \centering
    \includegraphics[scale=.5]{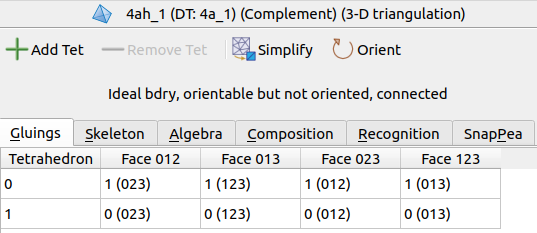}
    \caption{An ordered \textit{\textit{Regina}} triangulation of the complement of the figure-eight knot $4_1$}
    \label{fig:Regina:41:ordered}
\end{figure}

Figure \ref{fig:Regina:41:ordered} displays the same triangulation as Figure \ref{fig:Regina:41}, but with a different numbering of the vertices which makes it ordered. More precisely, indices $1$ and $2$ were swapped in the tetrahedron $0$, as can be seen by comparing the bottom rows of Figures \ref{fig:Regina:41} and \ref{fig:Regina:41:ordered}. 

In \textit{\textit{Regina}}, one can use to corresponding functions:
\begin{itemize}
\item \textbf{.isOrdered()} returns True if a given triangulation is ordered and False otherwise.
    \item \textbf{.order()} tries to find a permutation of indices which makes the triangulation ordered, returns False if the triangulation is not orderable, and returns True and a new modified ordered triangulation if it is orderable.
\end{itemize}

Figure \ref{fig:41:face:matrices}  represents the same ordered triangulation of $S^3\setminus 4_1$ where the faces were renamed with letters and additional information we will cover later. This figure illustrates an advantage of ordered triangulations: we can display them as such, just naming each face with a letter, without detailing the correspondence of triplets under gluing, since it is imposed by the fact that $X$ is ordered.

\begin{remark}
    \label{rmk:order-involution}
    \begin{enumerate}
        \item Note that an order is equivalent to an orientation of each quotient edge such that no face is a cycle. Indeed, given an order, each edge is oriented from the lower index vertex toward the upper index vertex; this is consistent with gluings and does not introduce cycles. Conversely, given such an orientation, in each tetrahedron, exactly one vertex is the source of all three edges, one the source of exactly one edge, one the source of exactly two edges and one the source of no edge. Label these vertices respectively $0,1,2,3$ and this gives an order.
        Such an order is also called a \textit{branching} in the literature (see for instance \cite{BB}).
        
        \item From the previous point, one sees that one can \emph{reverse} an order: in each tetrahedron, one can change the labeling of the vertices $0$, $1$, $2$, $3$ to $3$, $2$, $1$, $0$. Equivalently, one can reverse all (quotient) edges and this will not create a cycle.
    \end{enumerate}
\end{remark}

\subsection{Tetrahedron signs}

In Figure \ref{fig:41:face:matrices}, both tetrahedra are drawn as embedded in $\R^3$ (the $0$ vertex is the tip of the pyramid seen from above). Up to rotation in $\R^3$, there are two possible such embeddings, whether the triplet of vectors $(\overrightarrow{01}, \overrightarrow{02}, \overrightarrow{03})$ follows the right-hand rule or not. If it does follow the right-hand rule, we will say the tetrahedron $T$ has sign $\epsilon(T)=1$, and that $T$ has sign $\epsilon(T)=-1$ otherwise. We can then define a diagonal matrix
$\mathcal{E}=Diag(\epsilon(T_1),\ldots,\epsilon(T_N))$ associated to the spatial realization of $X$ (see Figure \ref{fig:41:face:matrices}). Because of the rigidity of face gluings, fixing the sign of one tetrahedron fixes every other sign, and so $\mathcal{E}$ is defined up to multiplication by $\pm 1$, according to the chosen spatial realization.

\begin{remark}\label{rmk:sign-involution}
Note that if the order is reversed (see Remark \ref{rmk:order-involution}), then the sign of each tetrahedron remains the same.
\end{remark}

\subsection{Face adjacency matrices}\label{sub:face:adj}

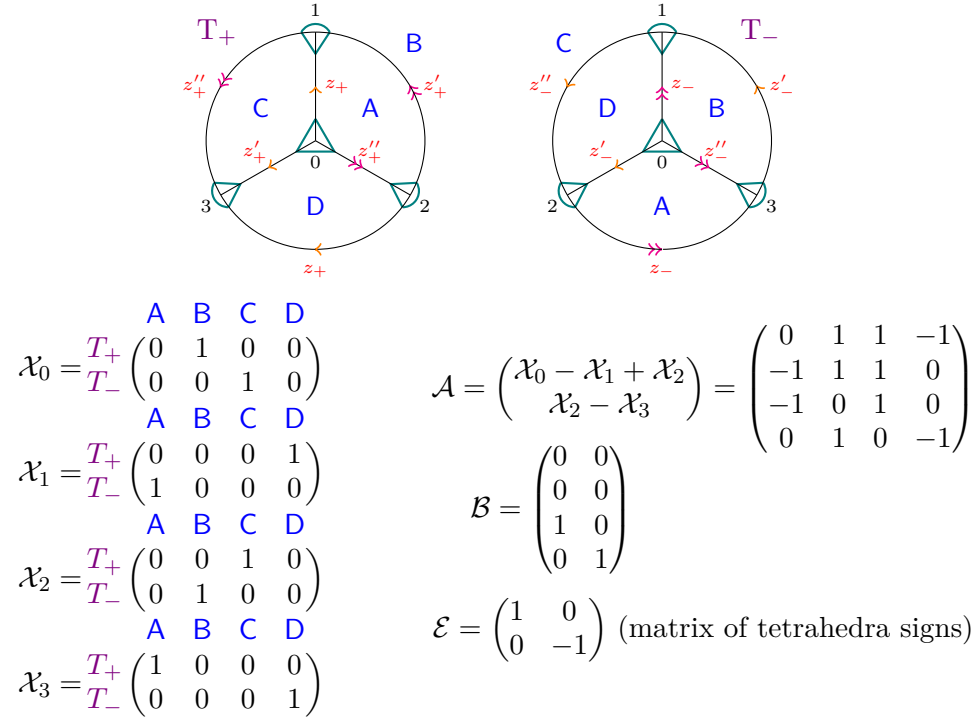
\begin{figure}[!h]
    \centering
\begin{tikzpicture}

\begin{scope}[scale=1.2]

\begin{scope}[xshift=0cm,yshift=0cm,rotate=0,scale=1.2]

\draw[color=violet] (-0.9,1) node{T$_+$} ;

\draw[color=red] (0.2,0.5) node{\tiny $z_+$} ;
\draw[color=red] (0,-1.2) node{\tiny $z_+$} ;
\draw[color=red] (0.5,-0.1) node{\tiny $z''_+$} ;
\draw[color=red] (-1.1,0.5) node{\tiny $z''_+$} ;
\draw[color=red] (-0.55,-0.1) node{\tiny $z'_+$} ;
\draw[color=red] (1.1,0.5) node{\tiny $z'_+$} ;

\draw (0,-0.2) node{\tiny $0$} ;
\draw (0,1.2) node{\tiny $1$} ;
\draw (1,-0.6) node{\tiny $2$} ;
\draw (-1,-0.6) node{\tiny $3$} ;

\begin{scope}[xshift=0cm,yshift=0cm,rotate=0,scale=.1]
    \draw[color=teal,thick] (0,2)--(-1.732,-1)--(1.732,-1)--(0,2);
\end{scope}

\begin{scope}[xshift=0cm,yshift=1cm,rotate=180,scale=.1]
    \draw[color=teal,thick ] (-1.3,0)--(0,2)--(1.3,0);
    \draw[color=teal,thick] (-1.3,0) arc (-150:-30:1.5);
\end{scope}

\begin{scope}[xshift=-0.86cm,yshift=-.5cm,rotate=-60,scale=.1]
    \draw[color=teal,thick ] (-1.3,0)--(0,2)--(1.3,0);
    \draw[color=teal,thick] (-1.3,0) arc (-150:-30:1.5);
\end{scope}

\begin{scope}[xshift=0.86cm,yshift=-.5cm,rotate=60,scale=.1]
    \draw[color=teal,thick ] (-1.3,0)--(0,2)--(1.3,0);
    \draw[color=teal,thick] (-1.3,0) arc (-150:-30:1.5);
\end{scope}

\draw[color=blue] (0.9,0.9) node{\small $\mathsf{\color{blue}B}$} ;
\draw[color=blue] (0,-0.6) node{\small $\mathsf{\color{blue}D}$} ;
\draw[color=blue] (-0.5,0.3) node{\small $\mathsf{\color{blue}C}$} ;
\draw[color=blue] (0.5,0.3) node{\small $\mathsf{\color{blue}A}$} ;

\path [draw=black,postaction={on each segment={mid arrow =orange, thick}}]
(0,0)--(-1.732/2,-0.5);

\path [draw=black,postaction={on each segment={mid arrow =orange, thick}}]
(0,0)--(0,1);

\path [draw=black,postaction={on each segment={mid arrow d =magenta, thick}}]
(0,0)--(1.732/2,-0.5);

\draw(1.732/2,-0.5) arc (-30:-89.5:1);
\draw[->,color=orange, thick](0,-1) arc (-89.5:-90:1);
\draw (0,-1) arc (-90:-150:1);

\draw(0,1) arc (90:29.5:1);
\draw[->>,color=magenta, thick](1.732/2,0.5) arc (29.5:30:1);
\draw (1.732/2,0.5) arc (30:-30:1);

\draw(0,1) arc (90:150:1);
\draw[->>,color=magenta, thick](-1.732/2,0.5) arc (149.5:150:1);
\draw (-1.732/2,0.5) arc (150:210:1);

\end{scope}

\begin{scope}[xshift=3.8cm,yshift=0cm,rotate=0,scale=1.2]

\draw[color=violet] (0.9,1) node{T$_-$} ;

\draw[color=red] (0.2,0.5) node{\tiny $z_-$} ;
\draw[color=red] (0,-1.2) node{\tiny $z_-$} ;
\draw[color=red] (0.5,-0.1) node{\tiny $z''_-$} ;
\draw[color=red] (-1.1,0.5) node{\tiny $z''_-$} ;
\draw[color=red] (-0.55,-0.1) node{\tiny $z'_-$} ;
\draw[color=red] (1.1,0.5) node{\tiny $z'_-$} ;

\draw (0,-0.2) node{\tiny $0$} ;
\draw (0,1.2) node{\tiny $1$} ;
\draw (1,-0.6) node{\tiny $3$} ;
\draw (-1,-0.6) node{\tiny $2$} ;

\begin{scope}[xshift=0cm,yshift=0cm,rotate=0,scale=.1]
    \draw[color=teal,thick] (0,2)--(-1.732,-1)--(1.732,-1)--(0,2);
\end{scope}

\begin{scope}[xshift=0cm,yshift=1cm,rotate=180,scale=.1]
    \draw[color=teal,thick ] (-1.3,0)--(0,2)--(1.3,0);
    \draw[color=teal,thick] (-1.3,0) arc (-150:-30:1.5);
\end{scope}

\begin{scope}[xshift=-0.86cm,yshift=-.5cm,rotate=-60,scale=.1]
    \draw[color=teal,thick] (-1.3,0)--(0,2)--(1.3,0);
    \draw[color=teal,thick] (-1.3,0) arc (-150:-30:1.5);
\end{scope}

\begin{scope}[xshift=0.86cm,yshift=-.5cm,rotate=60,scale=.1]
    \draw[color=teal,thick ] (-1.3,0)--(0,2)--(1.3,0);
    \draw[color=teal,thick] (-1.3,0) arc (-150:-30:1.5);
\end{scope}

\draw[color=blue] (-0.9,0.9) node{\small $\mathsf{\color{blue}C}$} ;
\draw[color=blue] (0,-0.6) node{\small $\mathsf{\color{blue}A}$} ;
\draw[color=blue] (-0.5,0.3) node{\small $\mathsf{\color{blue}D}$} ;
\draw[color=blue] (0.5,0.3) node{\small $\mathsf{\color{blue}B}$} ;

\path [draw=black,postaction={on each segment={mid arrow =orange, thick}}]
(0,0)--(-1.732/2,-0.5);

\path [draw=black,postaction={on each segment={mid arrow d =magenta, thick}}]
(0,0)--(0,1);

\path [draw=black,postaction={on each segment={mid arrow d =magenta, thick}}]
(0,0)--(1.732/2,-0.5);

\draw(1.732/2,-0.5) arc (-30:-89.5:1);
\draw[<<-,color=magenta, thick] (0,-1) arc (-89.5:-90:1);
\draw(0,-1) arc (-90:-150:1);

\draw(0,1) arc (90:29.5:1);
\draw[->,color=orange, thick](1.732/2,0.5) arc (29.5:30:1);
\draw (1.732/2,0.5) arc (30:-30:1);

\draw(0,1) arc (90:149.5:1);
\draw[->,color=orange, thick](-1.732/2,0.5) arc (149.5:150:1);
\draw (-1.732/2,0.5) arc (150:210:1);

\end{scope}

\end{scope}

\end{tikzpicture}

\begin{tikzpicture}

\draw (0,0) node {
$\mathcal{X}_0=\begin{pNiceMatrix}[first-row, first-col]
 &  \mathsf{\color{blue}A} & \mathsf{\color{blue}B} & \mathsf{\color{blue}C} & \mathsf{\color{blue}D} \\
\textcolor{violet}{T_+} & 0 & 1 & 0 & 0 \\
\textcolor{violet}{T_-} & 0 & 0 & 1 & 0
\end{pNiceMatrix}$
};

\draw (0,-1.4) node {
$\mathcal{X}_1=\begin{pNiceMatrix}[first-row, first-col]
 &  \mathsf{\color{blue}A} & \mathsf{\color{blue}B} & \mathsf{\color{blue}C} & \mathsf{\color{blue}D} \\
\textcolor{violet}{T_+} & 0 & 0 & 0 & 1 \\
\textcolor{violet}{T_-} & 1 & 0 & 0 & 0
\end{pNiceMatrix}$};

\draw (0,-2.8) node {
$\mathcal{X}_2=\begin{pNiceMatrix}[first-row, first-col]
 &  \mathsf{\color{blue}A} & \mathsf{\color{blue}B} & \mathsf{\color{blue}C} & \mathsf{\color{blue}D} \\
\textcolor{violet}{T_+} & 0 & 0 & 1 & 0 \\
\textcolor{violet}{T_-} & 0 & 1 & 0 & 0
\end{pNiceMatrix}$};

\draw (0,-4.2) node {
$\mathcal{X}_3=\begin{pNiceMatrix}[first-row, first-col]
 &  \mathsf{\color{blue}A} & \mathsf{\color{blue}B} & \mathsf{\color{blue}C} & \mathsf{\color{blue}D} \\
\textcolor{violet}{T_+} & 1 & 0 & 0 & 0 \\
\textcolor{violet}{T_-} & 0 & 0 & 0 & 1
\end{pNiceMatrix}$};

\draw (7,-.5) node {
$\mathcal{A} = \begin{pmatrix}
    \mathcal{X}_0-\mathcal{X}_1+\mathcal{X}_2\\\mathcal{X}_2-\mathcal{X}_3
\end{pmatrix}= \begin{pmatrix}
0&1&1&-1\\
-1&1&1&0\\
-1&0&1&0\\
0&1&0&-1
\end{pmatrix}$};

\draw (5,-1.6-.5) node {
$\mathcal{B}= \begin{pmatrix}
0&0\\
0&0\\
1&0\\
0&1
\end{pmatrix}$};

\draw (7,-3.2-.5) node {
$\mathcal{E}=\begin{pmatrix}
    1&0\\0&-1
\end{pmatrix}$ (matrix of {tetrahedra signs})};

\end{tikzpicture}

\caption{Thurston's triangulation of $M=S^3 \setminus 4_1$, and the face adjacency matrices}
    \label{fig:41:face:matrices}
\end{figure}

In this section we assume $X$ is ordered.
For $k=0,1,2,3$, let $x_k : X^3 \to X^2$ be the map defined by sending a tetrahedron to its face that is opposite to the $k$-th vertex, and $\mathcal{X}_k \in M_{N,2N}(\Z)$ the matrix of coefficients
$(\mathcal{X}_k)_{i,j}:=\delta_{j\text{-th face},x_k(T_i)}$. Finally, define \textcolor{black}{matrices $\mathcal{B}:=\begin{pmatrix}
    0_N \\{\rm Id}_N
\end{pmatrix} \in M_{2N,N}(\Z)$
and $\mathcal{A} = \begin{pmatrix}
    \mathcal{X}_0-\mathcal{X}_1+\mathcal{X}_2\\\mathcal{X}_2-\mathcal{X}_3
\end{pmatrix} \in M_{2N,2N}(\Z)$}. We informally call these matrices "face adjacency matrices" associated to the ordered triangulation $X$. See Figure \ref{fig:41:face:matrices} for an example.

\subsection{Neumann-Zagier matrices}
In this section, let us denote
$X^1=\{E_1,\dots, E_N\}$ is the set of 1-cells (edges) of the ideal triangulation $X$.

\subsubsection{Hyperbolic ideal tetrahedra}

Consider the upper half-space model 
of the hyperbolic 3-space $\mathbb{H}^3 = \C \times \R_{>0}$, whose boundary is $\C \cup \{\infty\}$. An ideal geodesic tetrahedron $T \subset \mathbb{H}^3$ can be assumed, up to isometry, to have ideal vertices $0,1,\infty$ and $z \in \C$;
we assume that $z \notin \{0,1\}$ and
we call $z$ the \textit{complex shape} of $T$. Up to cyclic permutation of the vertices, the same tetrahedron $T$ can be represented by the shapes $z'=\frac{1}{1-z}$ or $z''=\frac{z-1}{z}$. Recall or observe that $z z' z'' = -1$.

Equivalently, $T$ can be described by \textit{logarithmic shape parameters} 
$$\Log(z),\Log(z'),\Log(z'') \in \R +i(-\pi,\pi],$$ whose sum is $i\pi$. 

Finally, by writing $a:=\arg(z), b:=\arg(z'), c:=\arg(z'')$, we can describe $T$ by the triplet of real angles $(a,b,c)\in[0,\pi]^3$ of sum $\pi$. The angles $a,b,c$ are the values of the dihedral angles of $T$, where each pair of opposite edges has the same dihedral angle.

Assuming a shape structure is assigned to each tetrahedron $T_i$ of $X$,
we represent the complex shapes $z_i,z'_i,z''_i$ on associated edges of the triangulation (see an example in Figure \ref{fig:41:face:matrices}). If $X$ is ordered, we put $z_i$ on the edges $\overrightarrow{01}$ and $\overrightarrow{23}$ of $T_i$, and we require that $z_i,z'_i,z''_i$ have to cycle positively around each vertex in this order.
We denote $\textbf{z},\textbf{z'},\textbf{z''},\textbf{Log}\mathbf{(z)} ,\textbf{Log}\mathbf{(z')} ,\textbf{Log}\mathbf{(z'')} $ the associated vectors of $N$ coordinates.

\subsubsection{Weight functions, edge equations and angle structures}

For any edge $E_j \in X^1$, we define the associated \textit{complex weight function} 
$$\omega^\C_{X,E_j}(\textbf{Log}\mathbf{(z)} ,\textbf{Log}\mathbf{(z')} ,\textbf{Log}\mathbf{(z'')} )$$ as the sum of the logarithmic shapes $\Log(z_k^{\star})$ associated to all edges $e$ of the pre-gluing triangulation such that $[e]_\sim=E_j$.
See Figures \ref{fig:41:face:matrices} and \ref{fig:41:cusp:triang:NZ} for two examples of complex weight functions
$\omega^\C_{X,\textcolor{orange}{\rightarrow}}(\mathbf{Log(z^\cdot)})$,
$\omega^\C_{X,\textcolor{magenta}{\twoheadrightarrow}}(\mathbf{Log(z^\cdot)})$
for the edges $\textcolor{orange}{\rightarrow},\textcolor{magenta}{\twoheadrightarrow}$ in the triangulation of $S^3\setminus 4_1$.

Similarly, we can associate to any edge $E_j$ its \textit{angular weight function} $\omega^\R_{X,E_j}(\alpha)$, where $\alpha=(a_1,b_1,c_1,\ldots,a_N,b_N,c_N) \in [0,\pi]^{3N}$, as the sum of dihedral angles associated to all edges $e$ of the pre-gluing triangulation such that $[e]_\sim=E_j$. Observe that if $\textbf{Log(z)},\textbf{Log(z')},\textbf{Log(z'')}$ are the logarithmic shape parameters associated to $\alpha$, one has
$$\omega^\R_{X,E_j}(\alpha)= \Im \left(\omega^\C_{X,E_j}(\textbf{Log}\mathbf{(z)} ,\textbf{Log}\mathbf{(z')} ,\textbf{Log}\mathbf{(z'')} )\right).$$
The complex (resp. angular) \textit{edge equation} associated to $E_j$ is
$$\omega^\C_{X,E_j}(\textbf{Log}\mathbf{(z)} ,\textbf{Log}\mathbf{(z')} ,\textbf{Log}\mathbf{(z'')} )=2i\pi$$
(resp. $\omega^\R_{X,E_j}(\alpha)=2\pi$).

The \textit{space of angle structures} $\mathscr{A}_X$  associated to $X$ (which should not be confused with the face adjacency matrix $\mathcal{A}$ of Section \ref{sub:face:adj}) is defined as
$$\mathscr{A}_X:= \{\alpha=(a_1,\ldots,c_N) \in \fathi{(0,\pi)}^{3N} \ ; \
\forall j \in \{1,\ldots,N\}, \
\omega^\R_{X,E_j}(\alpha)=2\pi \}.
$$

\subsubsection{Cusp triangulation and holonomies}

By truncating the lone (post-gluing) vertex of $X$, we obtain a triangulation of $\overline{M}$ by truncated tetrahedra. This induces a triangulation of the boundary torus $\partial M$ called the \textit{cusp triangulation}, as can be seen in Figures \ref{fig:41:face:matrices} and \ref{fig:41:cusp:triang:NZ} (where the truncating triangles are drawn in teal).

\begin{figure}[!h]
    \centering

\begin{tikzpicture}

\begin{scope}[scale=1.2]

\draw (0,0)--(8,0)--(9,1.732)--(1,1.732)--(0,0);
\draw (1,1.732)--(2,0)--(3,1.732)--(4,0)--(5,1.732)--(6,0)--(7,1.732)--(8,0);

\filldraw[color=magenta] (0,0) circle (0.08cm);
\filldraw[color=magenta] (1,1.732) circle (0.08cm);
\filldraw[color=magenta] (0+8,0) circle (0.08cm);
\filldraw[color=magenta] (1+8,1.732) circle (0.08cm);
\draw[color=magenta,very thick] (0,0) circle (0.15cm);
\draw[color=magenta, very thick] (1,1.732) circle (0.15cm);
\draw[color=magenta,very thick] (8,0) circle (0.15cm);
\draw[color=magenta, very thick] (9,1.732) circle (0.15cm);
\filldraw[color=orange] (6,0) circle (0.08cm);
\filldraw[color=orange] (7,1.732) circle (0.08cm);
\draw[color=orange,very thick] (2,0) circle (0.08cm);
\draw[color=orange, very thick] (3,1.732) circle (0.08cm);
\draw[color=magenta,very thick] (4,0) circle (0.08cm);
\draw[color=magenta, very thick] (5,1.732) circle (0.08cm);
\draw[color=magenta,very thick] (4,0) circle (0.15cm);
\draw[color=magenta, very thick] (5,1.732) circle (0.15cm);

\draw[color=teal] (1,.6) node{$2_+$};
\draw[color=teal] (3,.6) node{$0_+$};
\draw[color=teal] (5,.6) node{$1_+$};
\draw[color=teal] (7,.6) node{$3_+$};
\draw[color=teal] (2,1.2) node{$1_-$};
\draw[color=teal] (4,1.2) node{$0_-$};
\draw[color=teal] (6,1.2) node{$2_-$};
\draw[color=teal] (8,1.2) node{$3_-$};

\draw[color=blue] (1,0) node{$\mathsf{\color{blue}D}$};
\draw[color=blue] (3,0) node{$\mathsf{\color{blue}D}$};
\draw[color=blue] (5,0) node{$\mathsf{\color{blue}C}$};
\draw[color=blue] (7,0) node{$\mathsf{\color{blue}C}$};
\draw[color=blue] (.5,.9) node{$\mathsf{\color{blue}A}$};
\draw[color=blue] (1.5,.9) node{$\mathsf{\color{blue}B}$};
\draw[color=blue] (2.5,.9) node{$\mathsf{\color{blue}C}$};
\draw[color=blue] (3.5,.9) node{$\mathsf{\color{blue}A}$};
\draw[color=blue] (4.5,.9) node{$\mathsf{\color{blue}B}$};
\draw[color=blue] (5.5,.9) node{$\mathsf{\color{blue}A}$};
\draw[color=blue] (6.5,.9) node{$\mathsf{\color{blue}D}$};
\draw[color=blue] (7.5,.9) node{$\mathsf{\color{blue}B}$};
\draw[color=blue] (8.5,.9) node{$\mathsf{\color{blue}A}$};
\draw[color=blue] (2,1.732) node{$\mathsf{\color{blue}D}$};
\draw[color=blue] (4,1.732) node{$\mathsf{\color{blue}D}$};
\draw[color=blue] (6,1.732) node{$\mathsf{\color{blue}C}$};
\draw[color=blue] (8,1.732) node{$\mathsf{\color{blue}C}$};

\draw[color=red] (.35,.2) node{\scriptsize $z''_+$};
\draw[color=red] (2.35,.2) node{\scriptsize $z'_+$};
\draw[color=red] (4.35,.2) node{\scriptsize $z''_+$};
\draw[color=red] (6.35,.2) node{\scriptsize $z'_+$};
\draw[color=red] (2-.35,.2) node{\scriptsize $z_+$};
\draw[color=red] (4-.35,.2) node{\scriptsize $z''_+$};
\draw[color=red] (6-.35,.2) node{\scriptsize $z_+$};
\draw[color=red] (8-.35,.2) node{\scriptsize $z''_+$};
\draw[color=red] (2,.45) node{\scriptsize $z'_-$};
\draw[color=red] (4,.45) node{\scriptsize $z''_-$};
\draw[color=red] (6,.45) node{\scriptsize $z'_-$};
\draw[color=red] (8,.45) node{\scriptsize $z''_-$};
\draw[color=red] (1,1.732-.45) node{\scriptsize $z'_+$};
\draw[color=red] (3,1.732-.45) node{\scriptsize $z_+$};
\draw[color=red] (5,1.732-.45) node{\scriptsize $z'_+$};
\draw[color=red] (7,1.732-.45) node{\scriptsize $z_+$};
\draw[color=red] (1+.35,1.732-.2) node{\scriptsize $z_-$};
\draw[color=red] (3+.35,1.732-.2) node{\scriptsize $z'_-$};
\draw[color=red] (5+.35,1.732-.2) node{\scriptsize $z_-$};
\draw[color=red] (7+.35,1.732-.2) node{\scriptsize $z'_-$};
\draw[color=red] (3-.35,1.732-.2) node{\scriptsize $z''_-$};
\draw[color=red] (5-.35,1.732-.2) node{\scriptsize $z_-$};
\draw[color=red] (7-.35,1.732-.2) node{\scriptsize $z''_-$};
\draw[color=red] (9-.35,1.732-.2) node{\scriptsize $z_-$};

\begin{scope}[xshift=4.2cm,yshift=2.3cm,rotate=-120]
\draw[color=olive,->] (0,0)--(3.2,0);    
\draw[color=olive] (3.5,0) node {$y$};
\end{scope}

\begin{scope}[xshift=0cm,yshift=1.075cm,rotate=0]
\draw[color=olive,->] (0,0)--(9.5,0);    
\draw[color=olive] (9.7,0) node {$x$};
\end{scope}

\end{scope}

\end{tikzpicture}

$(\textcolor{orange}{\rightarrow})$
starts at
\begin{tikzpicture}
\draw[color=orange, very thick] (0,0) circle (0.08cm);
\end{tikzpicture}
and ends at
\begin{tikzpicture}
\filldraw[color=orange] (0,0) circle (0.08cm);
\end{tikzpicture}
\ :\\ \color{red}
$
\omega^\C_{X,\textcolor{orange}{\rightarrow}}(\mathbf{Log(z^\cdot)})=
2 {\rm Log}(z_+)+{\rm Log}(z'_+)+2{\rm Log}(z'_-)+{\rm Log}(z''_-)=2i\pi$
\color{black}

$(\textcolor{magenta}{\twoheadrightarrow})$
starts at
\begin{tikzpicture}
\draw[color=magenta,very thick] (0,0) circle (0.08cm);
\draw[color=magenta,very thick] (0,0) circle (0.15cm);
\end{tikzpicture}
and ends at
\begin{tikzpicture}
\filldraw[color=magenta] (0,0) circle (0.08cm);
\draw[color=magenta,very thick] (0,0) circle (0.15cm);
\end{tikzpicture}
\ : \\
\color{red}
$ \omega^\C_{X,\textcolor{magenta}{\twoheadrightarrow}}(\mathbf{Log(z^\cdot)})=
{\rm Log}(z'_+)+2{\rm Log}(z''_+)+2{\rm Log}(z_-)+{\rm Log}(z''_-)=2i\pi$
\color{black}

\vspace{.2cm}

\textcolor{olive}{Meridian $y$} : \ 
\color{red}
$H^\C_{X,y}(\mathbf{z})= -{\rm Log}(z'_-)+{\rm Log}(z''_+)$
\color{black}

\textcolor{olive}{\textbf{Preferred longitude} $x+2y$} : \ 
\color{red}
$ H^\C_{X,x+2y}(\mathbf{z})=2i\pi -4{\rm Log}(z'_-)-2{\rm Log}(z''_-)$
\color{black}

\ 

$\mathbf{G}=\begin{pNiceMatrix}[first-row, first-col]
 &  {\color{red}z_+} & {\color{red}z_-} \\
{(\color{orange}\to\color{black})} & 2 & 0  \\
{\color{olive}x+2y} & 0 & 0 
\end{pNiceMatrix}, \ \ \
\mathbf{G'}=\begin{pNiceMatrix}[first-row, first-col]
 &  {\color{red}z'_+} & {\color{red}z'_-} \\
{(\color{orange}\to\color{black})} & 1 & 2  \\
{\color{olive}x+2y} & 0 & -4 
\end{pNiceMatrix}, \ \ \ 
\mathbf{G''}=\begin{pNiceMatrix}[first-row, first-col]
 &  {\color{red}z''_+} & {\color{red}z''_-} \\
{(\color{orange}\to\color{black})} & 0 & 1  \\
{\color{olive}x+2y} & 0 & -2 
\end{pNiceMatrix}$

\

$\mathbf{A}=\mathbf{G}-\mathbf{G'}=\begin{pNiceMatrix}[first-row, first-col]
 &  {\color{red}z_+} & {\color{red}z_-} \\
{(\color{orange}\to\color{black})} & 1 & -2  \\
{\color{olive}x+2y} & 0 & 4 
\end{pNiceMatrix}, \ \ \
\mathbf{B}=\mathbf{G''}-\mathbf{G'}=\begin{pNiceMatrix}[first-row, first-col]
 &  {\color{red}z''_+} & {\color{red}z''_-} \\
{(\color{orange}\to\color{black})} & -1 & -1  \\
{\color{olive}x+2y} & 0 & 2 
\end{pNiceMatrix}$

    \caption{Cusp triangulation of $S^3 \setminus 4_1$, and the Neumann-Zagier matrices}
    \label{fig:41:cusp:triang:NZ}
\end{figure}

Up to isotopy, any simple closed curve $\gamma$ in $\partial M$ can be made transverse to the cusp triangulation, and the intersection of $\gamma$ with each triangle separates one of the three angles from the other two. We can then define the \textit{complex holonomy} $H^\C_{X,\gamma}(\mathbf{z})$ of $\gamma$ 
 as the sum of terms of the form $\pm \Log(z_i),\pm \Log(z'_i),\pm \Log(z''_i)$, each logarithmic shape parameter corresponding to the lone angle of the previous sentence, and $\pm$ is $1$ if the piece of $\gamma$ circles this angle positively (and $-1$ otherwise). See Figure \ref{fig:41:cusp:triang:NZ}, where the complex holonomies of the meridian curve $y$ and the preferred longitude $x+2y$ are written.
 
 As for edge weight functions and edge equations, for $\alpha \in \mathscr{A}_X$, we define 
the \textit{angular holonomy} $H^\R_{X,\gamma}(\alpha)$ of $\gamma$ as the imaginary part of 
$H^\C_{X,\gamma}(\mathbf{z})$, where $\mathbf{z}$ is associated to $\alpha$. 

\subsubsection{Neumann-Zagier matrices}

 It is known that for a $3$-manifold with one cusp such as a knot complement, there exists a set of $N-1$ independent edges, in the sense that if the edge equations of these edges are satisfied, the last edge equation will automatically be satisfied.

Now we fix a simple closed curve $l$ in the boundary torus $\partial M$ (for a knot complement this will be the preferred longitude of the knot, which is null-homologous in the knot complement).

Let $\xi \in\C$ be a complex number.
We consider the set of $N$ equations on logarithmic shape parameters 
$\textbf{Log}\mathbf{(z)} ,\textbf{Log}\mathbf{(z')} ,\textbf{Log}\mathbf{(z'')} $ 
(where $\mathbf{Log(z^{\star})}$ is the vector of parameters $\Log(z^{\star}_k)$) given by
\begin{itemize}
    \item the first $N-1$ edge equations coming from edges of $X^1$ and
    \item the holonomy equation $H^\C_{X,l}(\mathbf{z}) = \xi$.
\end{itemize}
This system of equations can be written
\begin{align}\label{equ}
 \mathbf{G} \mathbf{Log(z)} +\mathbf{G'} \mathbf{Log(z')} + \mathbf{G''} \mathbf{Log(z'')} =
\begin{pmatrix}
    2i\pi \\ \ldots \\ 2i\pi \\ 
    \xi 
\end{pmatrix},
\end{align}
where $\mathbf{G},\mathbf{G'},\mathbf{G''} \in M_N(\Z)$ are the \textit{Neumann-Zagier matrices} associated to the previous numbering of edges and the choice of curve $l$.

We say that $X$ is \textit{geometric} when the previous system admits a solution $\mathbf{z}$ for $\xi=0$, where all the shape parameters have positive imaginary part (which means that all tetrahedra have positive hyperbolic volume).

Furthermore, we define the matrices
$$
\mathbf{A}:=\mathbf{G}-\mathbf{G'} \ \text{and} \ 
\mathbf{B}:=\mathbf{G''}-\mathbf{G'},
$$
which comes from rewriting the system on variables $\Log(z_k),\Log(z''_k)$, using the relation
$\Log(z'_k)=i\pi-\Log(z_k)-\Log(z''_k)$.

Getting back to the example of $S^3 \setminus 4_1$, consider Figure \ref{fig:41:cusp:triang:NZ}.
Edge equations associated to the orange simple arrow and the pink double arrow can be read on the cusp triangulation by going around a vertex (which represents a truncated edge) and are detailed below the picture. The same can be done for the holonomy equations associated to the meridian $y$ or the preferred longitude $l=x+2y$, which is the chosen curve here. The associated Neumann-Zagier matrices $\textcolor{black}{\mathbf{G},\mathbf{G'},\mathbf{G''},\mathbf{A},\mathbf{B}}$ are then detailed.

\subsection{FAMED triangulations}\label{sub:FAMED}

Let $l$ be a simple closed curve in the boundary torus $\partial M$.

We say that an ordered triangulation $X$ is FAMED (for "Face Adjacency Matrices with Edge Duality") with respect to the curve $l$ if 
\begin{enumerate}
\item the space of angle structures  is non empty: $\mathscr{A}_X \neq \emptyset$,
\item $\det \mathcal{A} \neq 0$, 
\item $\det \mathbf{B}\neq 0$, 
\item $\mathbf{B}^{-1}\mathbf{A} = \mathcal{X}_0 \mathcal{A}^{-1} \mathcal{B} \mathcal{E} + \left ( \mathcal{X}_0 \mathcal{A}^{-1} \mathcal{B} \mathcal{E}\right )^{\!\top} + \frac{\mathcal{E}+{\rm Id}_N}{2}$.
\end{enumerate}

For example, let us show that the ordered triangulation of $S^3\setminus 4_1$ displayed in Figures \ref{fig:41:face:matrices} and \ref{fig:41:cusp:triang:NZ} is FAMED. The 6-uple  $\left(\frac{\pi}{3},\ldots,\frac{\pi}{3}\right)$ is an angle structure (the complete one in fact, which illustrates that $X$ is geometric), so (1) is verified. (2) and (3) can be immediately checked, and (4) too, where the common value is $\begin{pmatrix}
    -1 & 0 \\ 0 & 2
\end{pmatrix}$.

{In what follows, "FAMED" will always mean "FAMED for the preferred longitude $l$ of the knot".}

\section{Experimental evidences for the volume conjecture}\label{sec:coding}

We explain here how we can prove, using automatic computations, that more than 42~000 knot complements admit an ordered triangulation that is FAMED and geometric; thus proving the {Andersen--Kashaev} volume conjecture for this number of knot complements. The gist of the computation is that checking the FAMED condition for an ordered geometric triangulation is an easy combinatorial task. The problem lies in finding an order - if it exists - that turns the triangulation into a FAMED one; and retriangulating if no such order exists.

The step by step explanation of the experimentation is based on a SageMath module hosted in the repository \cite{FAMEDexploration} that relies on \textit{Regina} and \textit{Snappy}. The complete sample exploration is illustrated in the notebook \verb|illustration.ipynb| that can be found inside the \verb|examples/| directory of the repository.

\subsection{Description of the census}

The census we work with is the Hoste-Thistlethwaite census of knots and links up to 14 crossings, as included in \textit{Snappy}. We restrict ourselves to knot complements. Each triangulation in this census has a name of the type K3a1 or K8n1: K denotes that this is a knot complement, and after is its DT name: the first number is the number of crossings, a/n describes if it is alternating or not and the last number is an index.

The best indicator of the complexity of checking the FAMED condition is the number of tetrahedra of the triangulation. We focus here our attention on the triangulations in the census that contain up to 23 tetrahedra. For reasons explained below, going further would require - at this stage - extremely long computations.
With this condition, we are left with 42.136 triangulations, and the list contains all triangulations of knot complements up to 12 crossings. Numerous triangulations with 13 and 14 crossings also appear.

We take as a small but significant sample of what can happen the first 14 elements of the census, that is, all the knots up to 7 crossings as listed in Table \ref{table:sample}. In this Table, we also give the more classical label of a knot complement, as found using KnotInfo \cite{knotinfo}.

\begin{table}[ht]
    \centering
    \begin{tabular}{|c|c|c|}
    \toprule
        {\bf Name in HTlinks} & {\bf Name} &{\bf Information} \\
    \midrule
        K3a1 & $3_1$ & Trefoil knot complement, not hyperbolic\\
        K4a1 & $4_1$ & $8$-knot complement, hyperbolic\\
        K5a1 & $5_2$ & hyperbolic\\
        K5a2 & $5_1$ & Torus knot, not hyperbolic\\
        K6a1 & $6_3$ & hyperbolic \\
        K6a2 & $6_2$ & hyperbolic \\
        K6a3 & $6_1$ & hyperbolic \\
        K7a1 & $7_7$ & hyperbolic \\
        K7a2 & $7_6$ & hyperbolic \\
        K7a3 & $7_5$ & hyperbolic \\
        K7a4 & $7_2$ & hyperbolic \\
        K7a5 & $7_3$ & hyperbolic \\
        K7a6 & $7_4$ & hyperbolic \\
        K7a7 & $7_1$ & Torus knot, not hyperbolic \\
    \bottomrule
    \end{tabular}

\ 
    
    \caption{Our sample of knot complements}
    \label{table:sample}
\end{table}

\subsection{Hyperbolicity}

The first point to check is whether those manifolds are hyperbolic. \textit{Snappy} has a built-in \verb|verify_hyperbolicity| method that we use for this step. This function actually verifies that the triangulation given by \textit{Snappy} is geometric, see also \cite{HIKMOT}.

Among the sample, we find that three triangulations are not geometric: K3a1, K5a2, K7a7. They indeed correspond to non hyperbolic manifolds, respectively the complements of the torus knots $3_1, 5_1$ and $7_1$. Among the whole census
one can verify that all manifolds but 13 are hyperbolic and that the \textit{Snappy} triangulations of the hyperbolic ones are indeed geometric \cite{HIKMOT}. This leaves us with 42~123 \emph{hyperbolic} manifolds - i.e. the 42~136 mentioned before minus the 13 non hyperbolic ones - for which we try and prove that a FAMED triangulation exists.

\subsection{Ordering triangulations}

For a given triangulation, we need to order it. As mentioned above, \textit{Regina} has a built-in method to construct an order. But this specific order may not be FAMED. So we need to check all orders until we find one which is FAMED (hopefully). The FAMEDexploration module contains tools to construct all orders for a given triangulation.
Table \ref{tab:orderings} shows, for the triangulations in the sample, the number of orders and the number of FAMED orders.
Recall from Remark \ref{rmk:order-involution} that an order is equivalent to an orientation of each edge of the triangulation such that no face is a cycle. 

As an ideal triangulation of a knot complement with $\nu$ tetrahedra has $\nu$ edges, the number of orders is bounded above by $2^\nu$. Moreover, we experimentally observe that a given order is FAMED if and only if its reverse order (in the sense of Remark \ref{rmk:order-involution}) is FAMED. 
Thus, in computation, we arbitrarily fix the orientation of one edge to avoid doing
twice the work.

\begin{table}[ht]
\centering
\begin{tabular}{|l|c|c|c|}
\toprule
Name & 
\multicolumn{3}{c|}{Number of:} \\
& 
Tetrahedra & 
Orders & 
FAMED orders\\
\midrule
K3a1 & 2 & 2 & 0  \\
K4a1 & 2 & 4 & 4  \\
K5a1 & 3 & 4 & 4  \\
K5a2 & 3 & 0 & 0  \\
K6a1 & 6 & 14 & 8  \\
K6a2 & 5 & 6 & 6  \\
K6a3 & 4 & 6 & 6  \\
K7a1 & 8 & 32 & 0  \\
K7a2 & 8 & 20 & 8  \\
K7a3 & 7 & 16 & 2  \\
K7a4 & 4 & 4 & 4  \\
K7a5 & 5 & 4 & 4  \\
K7a6 & 6 & 8 & 8  \\
K7a7 & 4 & 2 & 0  \\
\bottomrule
\end{tabular}

\ 

\caption{Orders data for the sample of the first 14 triangulations.}\label{tab:orderings}
\end{table}

We see first that:
\begin{itemize}
\item The triangulations of K3a1, K7a7 admit orders but they cannot be FAMED as the
associated manifolds are not hyperbolic. Indeed, the FAMED condition requires that the space of angle structures is non empty, which implies hyperbolicity by a theorem of Casson, see \cite[Corollary 4.6]{Lackenby-WordHypDehSurgery} or \cite[Theorem 1.1]{FuterGueritaud}
\item Not all triangulations admit orders: for example, K5a2 - which is not geometric -  admits no orders. We can check on the census that all \emph{geometric} triangulations that appear do admit orders.
\end{itemize}

For the other cases, we would hope that the ordered triangulation is FAMED. Indeed, in most cases in our sample, all orders give a FAMED triangulation. The exceptions are the following:
\begin{itemize}
  \item For K6a1, K7a2, K7a3, some orders are FAMED, and some other are not. Those examples are a very good place to explore  the FAMED condition further and its dependence on the order.
  \item For K7a1, we have plenty of orders (32 to be exact), but none of them is FAMED.
\end{itemize}

\subsection{Re-triangulation}

Hoping that all hyperbolic knot complements admit a FAMED triangulation, we resort to our last trick: we try to retriangulate the same manifold, by applying Pachner moves, until we find a triangulation admitting a FAMED order.

\textit{Regina} offers a method to explore the graph of retriangulation, under an upper bound on the maximal number of tetrahedra of a retriangulation - the \verb|depth| parameter of \textit{Regina}'s \verb|retriangulate|. When retriangulating K7a1 at depth 2, one gets a triangulation that admits a FAMED order.

Hence, all hyperbolic manifolds in our sample of 14 triangulations admit a FAMED triangulation.

This method, applied to all manifolds with at most 23 tetrahedra, gives the following experimental theorem:
\begin{theorem}
  All 42~123 but 6 hyperbolic manifolds in Hoste-Thistlethwaite census of knot complements with at most 14 crossings and 23 tetrahedra admit a triangulation with an order which is FAMED, and certified geometric. This list comprises all knot complements with at most 12 crossings - 
  
  For all those triangulations of knot complements, the Andersen-Kashaev volume conjecture holds.
\end{theorem}
The full list of triangulations proving this experimental theorem is available at \href{https://perso.imj-prg.fr/antonin-guilloux/famed-exploration-of-andersen-kashaev-volume-conjecture/}{Guilloux's webpage} 
 as a Regina container. The missing 6 manifolds have 13 or 14 tetrahedra in their initial triangulation {and are} K13a1017, K14a3125, K14a4490, K14n1207, K14n10807, K14n11640. For those 6 manifolds, computation is still {ongoing but is very long}.

One note on computation time: exploring the graph of retriangulation is a very long task whose complexity seems to grow doubly exponentially with the allowed depth. The whole census - but few examples - can be treated in a matter of a week, using parallelization - we used 20 parallel threads. For a dozen of examples with 23 tetrahedra, we need to allow depth 3 retriangulation as no FAMED triangulation can be found at depth 2; the computation, for each such example, lasts more than a week. For the 6 missing cases, it appears to be much longer. It looks like we are hitting a complexity wall at the threshold of 24 tetrahedra: most examples can be proved to be FAMED in a few days, but an exhaustive exploration would likely take weeks, if not months. What is lacking, according to our understanding, is a strategy (or at least a rationale) for retriangulating. This involves a better understanding of which patterns in a triangulation prevent FAMEDness, and how well-chosen Pachner moves may lead to a FAMED triangulation.

\section{Interesting patterns}\label{sec:patterns}

Upon inspection of all ordered triangulations, we see some patterns emerge that are not fully understood.
The inspection of the hyperbolic manifolds in the sample is performed in the companion notebook \verb|illustration.ipynb| of \cite{FAMEDexploration}.

The experimental evidence is univocal, as those patterns hold for all hundreds of thousands of triangulations inspected. So we state it as conjectures. Recall that the \emph{nullity} of a matrix is the codimension of its image.
\begin{conjecture}
  Let $T$ be an ordered and geometric triangulation of a knot complement. Then we have:
  \begin{enumerate}
    \item The nullity of $\bf{B}$ is twice the nullity of $\mathcal{A}$. 
    \item When they are invertible, then $\det(\mathcal{A}) = \pm 1$ and $\det(\bf{B})$ $= \pm 2$.
    \item When they are invertible, the triangulation $T$ is FAMED.
  \end{enumerate}
  In particular, conditions (2), (3) and (4) of the FAMED condition are equivalent.
\end{conjecture}
This conjecture hints at a deeper level of relations between the matrices $\mathcal A$ and $\bf{B}$.

\section{Future prospects}\label{sec:future}

This exploration of a sizeable census leads to further questions and experiments, and we still hope that those ideas could lead to unconditional proofs of the Andersen--Kashaev volume conjecture.

A first direction could be to use the data already computed to explore in-depth the FAMED condition:
\begin{itemize}
    \item We feel that low-complexity examples can guide us toward a better understanding of the impact of orders and retriangulations. Ideally, we should be able to choose how to retriangulate to arrive at a new triangulation that admits a FAMED order.
    \item Non FAMED examples are a good place to study the relations between the matrices $\mathcal A$ and $\bf{B}$ even when they are not invertible. Ideally, one could imagine a version of the condition for which we would need much fewer retriangulations.
    \item Our data set of FAMED/non FAMED ordered triangulations is huge - several hundreds of thousands of cases - and easy to extend further. We hope to extract some statistical study from this dataset. It may also be one of the places where machine learning and/or AI assisted exploration could lead to understanding hard-to-spot patterns, following the example of \cite{Lackenby}.
\end{itemize}

Another direction we are following, in collaboration with F. Guéritaud, is to extend the FAMED condition to cusped three manifolds that are not knot complements in $S^3$. A starting point would be once-punctured torus bundles over the circle: for this family, we have both a good theoretical understanding and efficient combinatorial ways of building numerous examples. A difficulty in this direction is that we have to find equivalents of the meridian and preferred longitude of the knot when the manifold is no longer a knot complement in $S^3$. Indeed, the construction of the Neumann-Zagier matrices depends on the choice of a curve $l$, which defines the last line.

\bibliographystyle{alpha}
\bibliography{bibli}

@misc{regina,
    author = {Benjamin A. Burton and Ryan Budney and William Pettersson and others},
    title = {Regina: Software for low-dimensional topology},
    howpublished = {{\tt http://\allowbreak regina-normal.\allowbreak github.\allowbreak io/}},
    year = {1999--2025}}

@misc{SnapPy,
     author={Culler, Marc and Dunfield, Nathan M. and Goerner,
     Matthias and Weeks, Jeffrey R.},
     title={Snap{P}y, a computer program for studying the geometry and topology of $3$-manifolds},
     howpublished={Available at \url{http://snappy.computop.org} (2025)},
}

@manual{sagemath,
  Key          = {SageMath},
  Author       = {{The Sage Developers}},
  Title        = {{S}ageMath, the {S}age {M}athematics {S}oftware {S}ystem ({V}ersion 10.3)},
  note         = {{\tt https://www.sagemath.org}},
  Year         = {2025},
}

@misc{FAMEDexploration,
  author = {Guilloux, A.},
  title = {A FAMED Exploration module},
  year = {2025},
  publisher = {PlmLab},
  journal = {GitLab repository},
  howpublished = {\url{https://plmlab.math.cnrs.fr/guilloux/FAMEDexploration}},
  commit = {d8eba440cd6dbd2ac8f84bd9d9998e1e2d400a4e}
}

@misc{HIKMOT,
      title={Verified computations for hyperbolic 3-manifolds}, 
      author={Neil Hoffman and Kazuhiro Ichihara and Masahide Kashiwagi and Hidetoshi Masai and Shin'ichi Oishi and Akitoshi Takayasu},
      year={2013},
      eprint={1310.3410},
      archivePrefix={arXiv},
      primaryClass={math.GT},
      url={https://arxiv.org/abs/1310.3410}, 
}

@article{Lackenby-WordHypDehSurgery,
 author = {Lackenby, Marc},
 title = {Word hyperbolic {Dehn} surgery},
 fjournal = {Inventiones Mathematicae},
 journal = {Invent. Math.},
 issn = {0020-9910},
 volume = {140},
 number = {2},
 pages = {243--282},
 year = {2000},
 language = {English},
 doi = {10.1007/s002220000047},
 zbMATH = {1463429},
 Zbl = {0947.57016}
}

@incollection{FuterGueritaud,
 author = {Futer, David and Gu{\'e}ritaud, Fran{\c{c}}ois},
 title = {From angled triangulations to hyperbolic structures},
 booktitle = {Interactions between hyperbolic geometry, quantum topology and number theory. Proceedings of a workshop, June 3--13, 2009 and a conference, June 15--19, 2009, Columbia University, New York, NY, USA},
 isbn = {978-0-8218-4960-6},
 pages = {159--182},
 year = {2011},
 publisher = {Providence, RI: American Mathematical Society (AMS)},
 language = {English},
 zbMATH = {5954470},
 Zbl = {1236.57002}
}

@article{Lackenby,
 author = {Davies, Alex and Juh{\'a}sz, Andr{\'a}s and Lackenby, Marc and Toma{\v{s}}ev, Nenad},
 title = {The signature and cusp geometry of hyperbolic knots},
 fjournal = {Geometry \& Topology},
 journal = {Geom. Topol.},
 issn = {1465-3060},
 volume = {28},
 number = {5},
 pages = {2313--2343},
 year = {2024},
 language = {English},
 doi = {10.2140/gt.2024.28.2313},
 zbMATH = {7927931},
 Zbl = {1555.57003}
}

@misc{knotinfo,
Author = {Livingston, Charles and Moore, Allison H.},
howpublished = {URL: \url{knotinfo.org}},
Month = {December},
Title = {KnotInfo: Table of Knot Invariants},
Year = {2025},
}

@article{AK,
 author = {Andersen, J{\o}rgen Ellegaard and Kashaev, Rinat},
 title = {A {TQFT} from quantum {Teichm{\"u}ller} theory},
 fjournal = {Communications in Mathematical Physics},
 journal = {Commun. Math. Phys.},
 issn = {0010-3616},
 volume = {330},
 number = {3},
 pages = {887--934},
 year = {2014},
 language = {English},
 doi = {10.1007/s00220-014-2073-2},
 zbMATH = {6324352},
 Zbl = {1305.57024}
}

@misc{BAW,
 author = {Ben Aribi, Fathi and Wong, Ka Ho},
 title = {The {Andersen}-{Kashaev} volume conjecture for {FAMED} geometric triangulations},
 year = {2024},
 howpublished = {Preprint, {arXiv}:2410.10776 [math.{GT}] (2024)},
 url = {https://arxiv.org/abs/2410.10776},
 arXiv = {arXiv:2410.10776}
}

@book{MY,
 author = {Murakami, Hitoshi and Yokota, Yoshiyuki},
 title = {Volume conjecture for knots},
 fseries = {SpringerBriefs in Mathematical Physics},
 series = {SpringerBriefs Math. Phys.},
 issn = {2197-1757},
 volume = {30},
 isbn = {978-981-13-1149-9; 978-981-13-1150-5},
 year = {2018},
 publisher = {Singapore: Springer},
 language = {English},
 doi = {10.1007/978-981-13-1150-5},
 zbMATH = {6909717},
 Zbl = {1410.57001}
}

@article{BB,
 author = {Baseilhac, Stephane and Benedetti, Riccardo},
 title = {Quantum hyperbolic geometry},
 fjournal = {Algebraic \& Geometric Topology},
 journal = {Algebr. Geom. Topol.},
 issn = {1472-2747},
 volume = {7},
 pages = {845--917},
 year = {2007},
 language = {English},
 doi = {10.2140/agt.2007.7.845},
 zbMATH = {5220898},
 Zbl = {1139.57008}
}

\end{document}